\magnification=\magstep1
\baselineskip=16pt
\lineskip=0pt
\lineskiplimit=0pt
\parskip=6pt
\vsize=7.65in


\def\Z   {{\bf Z}}

\def\R   {{\bf R}}

\def\P   {{\bf P}}

\def\F   {{\bf F}}


\def\det   {{\rm det}}
\def\Adj   {{\rm Adj}}


\def\Kbar      {\overline{\cal K}}

\def\Kahler    {K\"ahler\ }

\def\pi        {\partial /\partial z_i}

\def\ip {\, \raise 1.5pt \hbox{\underbar{\hskip 4.3pt }\vrule height 1.5pt
depth 1.4pt width .2pt}\, }

\def\ha        {\textstyle {1\over 2}\displaystyle}

\def\qu        {\textstyle {1\over 4}\displaystyle}
\def\si        {\textstyle {1\over 6}\displaystyle}


\def\calK      {{\cal K}}

\noindent \bf Some Geometry and Combinatorics for the $S$-invariant of Ternary Cubics.\rm
\bigskip \noindent \rm P.M.H. WILSON \medskip
\noindent Department of Pure Mathematics, University of Cambridge,\hfil\break 16
Wilberforce Road, Cambridge CB3 0WB, UK
\hfil\break email : pmhw@dpmms.cam.ac.uk \medskip
\noindent Max-Planck-Institut f\"ur Mathematik,\hfil\break
 Vivatsgasse 7, 53111 Bonn, Germany

\noindent \it Abstract. \rm In earlier papers 
[Wilson 04, Totaro 04], the $S$-invariant of a ternary cubic
$f$ was interpreted in terms of 
the curvature of related Riemannian and pseudo-Riemannian metrics --- this 
is clarified further in Section 1.  In
the case when
$f$ arises from the cubic form on the second cohomology of a 
smooth projective threefold with second
Betti number three, the value of the $S$-invariant is closely linked to the behaviour of this
curvature on the open cone consisting of K\"ahler classes.  
In this paper, we concentrate on  the cubic forms arising
from complete intersection threefolds in the product of three projective spaces, and  investigate
various conjectures of a combinatorial nature arising from their invariants.

\noindent \it Keywords:  \rm ternary cubics, invariant theory, curvature, combinatorial inequalities
\hfill\break
\noindent Mathematics Subject Classification 2000: Primary 15A72, 32J27, Secondary 14H52, 53A15

\noindent \bf INTRODUCTION.

\rm Given a real  form $f(x_1 , \ldots , x_m )$ of degree $d >2$, there is a 
pseudo-Riemannian metric, given by the matrix 
$ (g_{ij}) = -  ( \partial ^2 f / \partial x_i \partial x_j )/d(d-1)$, defined on 
the open subset of $\R ^m$ where the determinant $h = \det (g_{ij})$ is non-zero.  
This metric is referred to by Totaro as the \it Hessian metric\rm , and we  
study it further in the case when $f$ is a real ternary cubic.
Building on previous work in [Wilson 04, Totaro 04], we determine the full curvature 
tensor of this metric in terms of $h$ and the $S$-invariant of $f$ 
(Theorem 1.3).  

Motivated by the 
geometric background, as summarised below, we are led 
to consider cubic
forms associated  to complete intersection threefolds in the product of three projective 
spaces, and from this to study those cubic forms which arise as follows:  
We choose positive integers $d_1 , d_2 , d_3$ and $r\ge 0$ such that 
$d_1 + d_2 + d_3 = r +3$, and set 
$$ P = (x_1 H_1 + x_2 H_2 + x_3 H_3)^3 \prod _{j=1} ^r (a_j H_1 + b_j H_2 + c_j H_3 ),$$ 
with the $a_j$, $b_j$, $c_j$ are non-negative, and such that the cubic 
$F(x_1 , x_2 , x_3)$, defined by taking the coefficient of the term in $H_1 ^{d_1} H_2 ^ {d_2} H_3
^{d_3}$ in the above formal product $P$, is non-degenerate.  
Calculations from [Wilson 04] and the further discussion provided below suggest various 
conjectures concerning the invariants of such cubics.  In this paper, we shall 
concentrate mainly on Conjecture 2.1 that, regarding the $S$-invariant
 as a polynomial in the 
$a_j$, $b_j$ and $c_j$,
every coefficient is non-negative.
Extensive computer investigations are described in support of this conjecture.

In Section 3, we consider the cofactors $B_{pq}$ of the Hessian matrix of $F$ 
given by the matrix of second partial derivatives.
In the 
specific case under consideration, these are polynomials in $x_1, x_2, x_3$ 
and the $a_j$, $b_j$, $c_j$.  
We derive formulae for the coefficients of these polynomials, 
and deduce that these
coefficients are negative if 
$p=q$ and positive if $p \ne q$ (Theorem 3.1).  From this, we deduce that 
the Hessian determinant 
$H$ of $F$, that is the determinant of the Hessian matrix,  only has positive
coefficients.  This latter result represents a combinatorial version of
the Hodge index theorem.

In the final section, we return to a formula for $S$, given in Section 1, in 
terms of the cofactors $B_{pq}$ of the Hessian matrix. 
The fact that for the cubics $F$ being considered, we have formulae 
for the coefficients of monomials in the $B_{pq}$,  
enables us to produce a 
useful algorithm for determining the coefficient of a given monomial in $S$.  
We  run this algorithm for some critical cases, where we check that the 
conjectured positivity holds.\medskip

\noindent \bf 0.  THE GEOMETRIC BACKGROUND.

\rm In this preliminary section, the theory and calculations of this paper are 
set in their geometric context, and motivation is given for the conjectures 
appearing in Section 2.

For a compact \Kahler $n$-fold $X$, we can consider the level set in 
$H^{1,1} (X, \R ) \subset H^2(X, \R )$ defined by setting the degree $n$ form 
$D \mapsto D^n$ (given by cup-product) to be one.  The intersection of 
this level set with 
the \Kahler cone $\calK \subset H^{1,1} (X, \R )$ gives a manifold 
$\calK _1$ of dimension $h^{1,1} -1$, on which there is a natural Riemannian metric. 
The tangent space to $\calK _1$ at a point $D$ may be identified as 
$\{ L \in H^{1,1} 
:  D^{n-1} \cdot L =0 \}$,  and the Riemannian metric specified by 
$$ (L_1 , L_2 ) \mapsto - D^{n-2} \cdot L_1 \cdot L_2 .$$  
This is precisely the restriction to $\calK _1$ of the Hessian metric (as defined above) associated to the
degree $n$ form on $H^{1,1} (X, \R )$.
In [Wilson 04], we initiated
the study of this manifold and its curvature, motivated mainly by the implications that
any restrictions on this curvature might have concerning the existence and
classification of Calabi--Yau threefolds with a given differentiable structure.

In the cited paper, we  showed that if one assumed the existence of limit points in
complex moduli corresponding to a certain specified type of degeneration, then the
sectional curvatures of $\calK _1$ were bounded between $- {1\over 2} n (n+1)$ and $0$.
In the particular case of complex projective threefolds with second betti number 3 and 
$h^{2,0} =0$, we have a ternary cubic form $F$ on $H^2(X, \R )$, and an explicit formula
was produced for the curvature of 
the surface $\calK _1$, namely 
$$- \textstyle {9\over 4}\displaystyle + \qu 6^6 S F^2/H^2, $$ 
where $S$ denotes the
$S$-invariant of $F$ (see Section 1 below) and $H$ the Hessian determinant.  From this
one notes that if $S \ne 0$ and there exists a point $D$ on the boundary of the \Kahler
cone at which $H$ vanishes but
$F$ doesn't,  then the curvature is unbounded on $\calK _1$.  

In the case of a Calabi--Yau threefold,  
a rational such point $D$ can be seen to correspond to the
contraction of a surface on 
$X$ to a point [Wilson 92].  If $D$ lies in the interior of a codimension one
face  of the closure $\Kbar$ of the \Kahler cone (recalling from [Wilson 92] that away
from $F=0$, the boundary of $\Kbar$ is locally rational polyhedral), then 
in appropriate coordinates the cubic form may be written as 
$F = a x_1^3 + g(x_2 , x_3 )$, and in particular $S=0$.  
If however $D$ generates an extremal ray of $\Kbar$ (i.e. corresponds to a codimension 2 
face of $\Kbar$), we automatically have that $D$ is rational.  Moreover we may have 
that $S$
is non-zero, although using the classification of contractions from [Wilson 92], one
can show that in this case $S$ must be non-negative.  
There exist examples of  such Calabi--Yau threefolds with $S>0$, and hence with
the curvature of
$\calK _1$  unbounded above --- the simplest examples here are  provided by general
Weierstrass  fibrations over $\P ^1 \times \P^1$ or 
over the Hirzebruch surfaces $\F _1$ and $\F _2$.  In 
the non Calabi--Yau case, an even simpler example is provided by taking the cone (in
$\P ^4$) on a smooth quadric surface in $\P ^3$, and blowing up in the singular point; 
the first of the Calabi--Yau examples given above is closely related to this one.  In the
examples 
above, the curvature is in fact strictly positive on $\calK _1$, but there are a number
of examples of Calabi--Yau hypersurfaces in weighted projective 4-space, with second 
betti number three, where the curvature tends to infinity as one approaches some
extremal ray on the 
boundary, but with it being negative at other points of $\calK _1$.  It follows 
however from the above above discussion,
at least in the case of Calabi--Yau threefolds with second betti number 3 and 
$h^{2,0} =0$, that the curvature of the surface $\calK _1$ is bounded below, 
and in the case when it is not bounded above this lower bound will be $-9/4$, or in
other words the $S$-invariant of the ternary cubic is non-negative. 

The ideas in [Wilson 04] were motivated in the case of Calabi--Yau $n$-folds by a Mirror
Symmetry argument relating the curvature on $\calK _1$ to the curvature of the 
Weil--Petersson metric on the complex moduli space of the mirror.  Known results on 
the curvature of the Weil--Petersson metric 
in fact involve the Ricci curvature rather than sectional curvatures and only in 
general provide a lower bound; for
3-folds and 4-folds one can however construct from the Weil--Petersson metric and its
Ricci curvature an associated metric, the \it Hodge metric\rm , and there are then both
upper and lower bounds for the sectional curvatures of this metric [Lu 01, Lu \& Sun 
04]. Thus one should not perhaps be  surprised by the examples given above where the
curvature of $\calK_1$ is  positive --- this is expected in the mirror to correspond to
the Ricci curvature  of the Weil--Petersson metric being positive (in some 
neighbourhood of a large complex limit point).  One should however expect a lower bound
for the  Ricci curvature of the metric on $\calK _1$, and in the Calabi--Yau threefold
case  the author conjectures this to be $-({n\over 2 })^2(h^{1,1} -2) $.  In fact using 
calculations from ([O'Neill 83], p. 211) and the interpretation of the Hessian 
metric in terms of a warped product (Lemma 2.1 of [Totaro 04]), this
conjecture may be checked  
to be equivalent to the semi-Riemannian Hessian metric on $\calK$  having \it
non-negative
\rm Ricci curvature (compare also with the explicit formula produced in Theorem 1.3). 
This latter rather attractive conjecture lends itself to being verified by computer, 
and has been checked by the author to hold for all
the standard examples given in [Wilson 04], and also for certain Calabi--Yau threefolds
with rather larger values of
$b_2$ --- one such example corresponding to a hypersurface of degree 13 in
weighted projective  space $\P (1, 2,3,3,4)$, a Calabi--Yau threefold with $b_2 =5$.

Thus for Calabi--Yau threefolds with $b_2 =3$ and $h^{2,0} =0$, the above expectation
corresponds to $\calK _1$ having curvature bounded below by $-9/4$, or equivalently 
to the $S$-invariant of the ternary cubic form being non-negative.  This has been 
extensively checked against available lists of Calabi--Yau threefolds with $b_2 =3$. 
It should be noted that the non-negativity of $S$ is known to fail in general for 
complex projective threefolds with $b_2 =3$ and $h^{2,0} =0$ [Wilson 04].
However, in the case of complex projective threefolds admitting the 
specific type of degeneration described in [Wilson 04], the author expects 
the lower bound will again be $-9/4$ rather than the $-3$  as proved there --- for
higher  values of $b_2$, the lower bound of $-3$
on the the \it sectional \rm curvature can be achieved.  In the case for instance of 
abelian threefolds, the lower bound of $-3$
on the \it sectional \rm curvature is attained, although one checks easily that the
sharper lower bound of $-9 (h^{1,1} -2)/4$ holds for the Ricci curvature.

For the general case of complete intersection threefolds in the product of three
projective  spaces, one has that $b_2 =3$ and $h^{2,0} =0$, and one can
degenerate the  defining polynomials into products of polynomials on the three factors,
and each of these polynomials may be assumed to be the  
product of distinct linear forms.  The author expects (but it
will be  non-trivial to prove) that the general such degeneration will be of the type 
described in [Wilson 04], with the product of harmonic two forms being approximately  
harmonic.  The conjectures introduced in Section 2 will then be closely related to 
the conjecture that the curvature of $\calK _1$ for such threefolds is bounded between 
$-9/4$ and zero.

The case of complete intersection threefolds in the product of three projective 
spaces therefore represents an important test case for the above conjectures and
speculations.  We shall see that they lead to rather striking positivity conjectures of
a combinatorial nature, involving the classical invariants of ternary cubic forms, for
which extensive computational evidence will be presented.

\noindent \bf 1. THE $S$-INVARIANT AND CURVATURE.

\rm We consider a general non-degenerate ternary cubic form with real coefficients 
$$\eqalign{ f = & a_{300}x_1 ^3 + a_{030}x_2 ^3 + 
a_{003}x_3 ^3 + 3 a_{210}x_1 ^2 x_2 + 
3 a_{201}x_1 ^2 x_3 \cr & + 3 a_{120} x_1 x_2 ^2 + 3 a_{021}x_2 ^2 x_3 + 
3 a_{102} x_1 x_3 ^2 + 3 a_{012} x_2 x_3 ^2 + 6 a_{111}x_1 x_2 x_3 . \cr} $$
Associated to $f$, we have two basic invariants $S$ and $T$, one of degree 4 
in the coefficients and one of degree 6 [Aronhold 58, Sturmfels 93].  
The $S$-invariant is given
explicitly (see [Sturmfels 93], page 167) by an expression in the coefficients 
with 25 terms 
$$\eqalign{S = & a_{300} a_{120} a_{021} a_{003} - a_{300} a_{120} a_{012}^2 - 
a_{300}a_{111}a_{030}a_{003} + a_{300}a_{111}a_{021}a_{012} 
\cr & + \quad \ldots \quad + a_{201}a_{111}a_{102}a_{030} + 
a_{120}^2 a_{102}^2  - 2 a_{120}a_{111}^2 a_{102} + a_{111}^4 .\cr}$$

As indicated above, this invariant is closely associated with curvature.  We define
the \it index cone \rm in $\R ^3$ to consist of the  points at which $f$ is positive and
the indefinite metric defined by the matrix $f_{ij} = 
\partial ^2 f / \partial x_i \partial x_j$ is of signature $(1,2)$.  The 
restriction of $g_{ij} = - {1\over 6} f_{ij}$ to the level set $M$ given by 
$f=1$ in the index cone 
is then a Riemannian metric, whose curvature at any  point is given
by the formula 
$$- \textstyle {9\over 4}\displaystyle + \qu  S f^2/h^2, $$  where 
$h = \det (g_{ij}) = - H/6^3$, with $H$ denoting the Hessian determinant of $f$ 
([Wilson 04], Theorem 5.1).  
Strictly speaking, we do not need to 
include the $f^2$ in this formula, since by definition it has 
value one on the level set; however 
for any point in the index cone, the formula given provides the curvature at
the unique point of $M$ on the corresponding ray.  This formula was both 
extended to higher degrees and clarified further in [Totaro 04].  

Consider now
the pseudo-Riemannian metric
defined by the matrix $ g_{ij} = 
- \si f_{ij}$, on a suitable open subset of $\R ^3$.  In the case 
of cubics, Theorem 3.1 of [Totaro 04] reduces to the following statement: if $U$ is
an open subset of $\R ^3$ on which the Hessian $H$ is non-zero, and $M$ denotes
the level set in $U$ given by $f=1$, then the sectional curvature of $U$ on 
the tangent 2-plane to $M$ at a point is just $6^6 S f/H^2 = S f/h^2$.  This reproves 
the formula given above for the curvature of the restricted metric to $M$ 
and generalises in a natural way to forms $f$ of arbitrary degrees $>2$ 
([Totaro 04], (3.1)).
It should be noted here that for ternary cubics $f$,
 the Clebsch version $S(f)$ of the $S$-invariant 
(as used in Totaro's paper) is the Aronhold $S$-invariant (as used in this paper) 
multiplied by a factor $6^4$.

One point that I wish to emphasize in this section is that, once one knows 
the $S$-invariant and the Hessian determinant $H$, the whole curvature tensor of the
above  pseudo-Riemannian metric is given very simply by (1.3), thus extending in this
case  Theorem 3.1 from [Totaro 04].  

Throughout this paper, we shall denote by $B$ the adjoint matrix to 
$A = (f_{ij})$, with entries the cofactors of $A$.  We shall need the 
following identity, proved by classical invariant theory.

\proclaim Lemma 1.1.  $$ \ha \sum _{p, q} B_{pq} 
(\partial ^2 B_{ij}/ \partial x_p \partial x_q) = 
6^4 S x_i x_j.$$\par

\noindent \it Proof.\quad  \rm If we apply the Clebsch polarization operator 
$\sum y_i \, {\partial /\partial x_i}$ to $f$ twice, we obtain a mixed concomitant 
$ S^3 V^* \to S^2 V^* \otimes V^*$ (where $V$ denotes  the 3-dimensional real vector space),
 which in
coordinates may be written as 
$$ f \mapsto \sum _{i,j} y_i \, y_j \, {{\partial ^2 f }\over {\partial x_i \partial x_j}} .$$
Passing to the dual quadratic form (scaled by $H( {\bf x })$), we 
obtain a mixed concomitant 
$ S^3 V^* \to S^2 V^* \otimes S^2 V$, which in coordinates may be written as 
$$ f \mapsto \sum _{p,q} B_{pq} ({\bf x})\, \partial /\partial y_p \, \partial / \partial y_q .$$

Taking a convolution of two such concomitants, contracting out a factor $S^2 V \otimes S^2 V^* $, 
we obtain a concomitant 
$ S^3 V^* \to S^2 V^* \otimes S^2 V$, which in coordinates may be written as 
$$ f \mapsto \sum _{i,j} \biggl ( \sum _{p,q} B_{pq} ({\bf x}) \ {{\partial^2 B_{ij} ({\bf x})}
\over {\partial x_p \partial x_q }} \biggr )
\, \partial /\partial z_i \, \partial /
\partial z_j .$$

We can check easily on the Hesse cubic $x_1 ^3 + x_2 ^3 + x_3 ^3 + 6\lambda \, x_1 x_2 x_3$ that 
$$ {1\over 2} \sum _{p,q} B_{pq} ({\bf x}) \ {{\partial^2 B_{ij} ({\bf x})}
\over {\partial x_p \partial x_q }} = 6^4 S\, x_i x_j ,$$ 
where $S = \lambda (\lambda ^3 -1)$ in the $S$-invariant, and hence we deduce that the two 
concomitants $$ \sum _{i,j} 
{1\over 2} \biggl ( \sum _{p,q} B_{pq} ({\bf x}) \ {{\partial^2 B_{ij} ({\bf
x})}
\over {\partial x_p \partial x_q }} \biggr )
\, \partial /\partial z_i \ \partial /
\partial z_j \quad  \hbox{\rm and}\quad \sum _{i,j} 6^4 S\, x_i x_j 
\ \partial /\partial z_i \ \partial /
\partial z_j $$ are identical, since clearly they also transform in the same way 
under the operation
of scaling the coordinates.  Thus we deduce the result claimed.
\bigskip

\noindent \bf Remark 1.2. \quad\rm  If we now express the cofactors $B_{pq}$ in terms 
of the $f_{ij}$, and then operate on both sides of (1.1) by 
$\partial ^2 / \partial x_i \partial x_j$, we get formulae for $S$ analogous to 
those given on page 116 of [Aronhold 58].
From (1.1), it  follows immediately that, for any $i,j$,  
$$  \ha \sum _{p, q} (\partial ^2 B_{pq} / \partial x_i \partial x_j) 
(\partial ^2 B_{ij}/ \partial x_p \partial x_q) = 
6^4 (1 + \delta _{ij})\, S .$$ 
\bigskip
It is shown in [Totaro 04] that the 
curvature tensor of the pseudo-Riemannian metric defined above has components 
$$R_{ijkl} = - {1\over {144}} \sum _{p,q} g^{pq} (f_{jlp} f_{ikq} - 
f _{ilp} f_{jkq} ),$$ where $(g^{pq})$ denotes the inverse matrix to $(g_{ij})$.
Thus, for instance, if we let $h = \det (g_{ij}) = - H/6^3$, then 
$$ - 4 h R_{1212} = 6^{-4}  \sum _{p,q} B_{pq} (f_{11p} f_{22q} - f _{12p} f_{12q} ).$$

We now observe that $$(f_{11p} f_{22q} + f_{11q} f_{22p} - 
2 f _{12p} f_{12q}) = 
\partial ^2 (f_{11} f_{22} - f_{12}^2 )/ \partial x_p \partial x_q = 
\partial ^2 B_{33}/ \partial x_p \partial x_q ,$$ and so 
$$ \eqalign{ - 4 h \, 6^4 R_{1212} & =  \ha \sum _p B_{pp} (\partial ^2 B_{33}/ \partial x_p
\partial x_p) +
 \sum _{p < q} B_{pq} (\partial ^2 B_{33}/ \partial x_p \partial x_q) \cr 
& =  \ha \sum _{p, q} B_{pq} (\partial ^2 B_{33}/ \partial x_p \partial x_q).\cr }$$ 
Hence we deduce from (1.1) that $ -4 h R_{1212} = S x_3^2 $.

In Lemma 1.1, we can also take $(i,j) = (1,2)$.  Since 
$B_{12} = f_{13} f_{23} - f_{12} f_{33}$, for any given 
$(p,q)$ we have $$\partial ^2 B_{12}/ \partial x_p \partial x_q = 
f_{13p} f_{23q} + f_{13q} f_{23p} - f_{12p}f_{33q} - f_{12q} f_{33p} .$$
The formula for curvature then implies that 
$$ \ha 6^{-4} \sum _{p,q} B_{pq} (\partial ^2 B_{12}/ \partial x_p \partial x_q ) 
=  4 h R_{1323} ,$$ and so we deduce from (1.1) that 
$ S x_1 x_2 = 4 h R_{1323} = - 4 h R_{1332} = -4 h R_{3123}$. 

\proclaim Theorem 1.3.  All components of the curvature tensor of the 
Hessian metric on $U \subset \R ^3$, where $U$ is the open subset 
given by the non-vanishing of $H$, are
given simply in terms the invariant
$S$ and the Hessian of $f$, and are all of the form $\pm \qu S \, x_i x_j /h$ 
for appropriate $i,j$ and choice of sign. Moreover, given tangent vectors 
$ \xi = \sum \lambda _i \, \partial / \partial x_i$ and 
$\eta = \sum \mu _j \, \partial / \partial x_j$, the corresponding value of the curvature 
tensor satisfies $$ - 4 h \, R(\xi ,\eta ,\xi ,\eta ) = 
S \, (\lambda _1 \mu _2 x_3 + 
\lambda _2 \mu _3 x_1  + \lambda _3 \mu _1 x_2 - \lambda _2 \mu_1 x_3 -
\lambda _3 \mu _2 x_1 - \lambda _1 \mu _3 x_2 )^2 .$$\par
\noindent \it Proof.\rm\quad  Since we have formulae for $R_{1212}$ and 
$R_{1323}$, we have the analogous formulae for $R_{ijij}$ and $R_{ijkj}$.  
We now use the general fact that the curvature tensor is invariant under 
exchanging the first pair of indices with the second pair of indices, and 
is anti-invariant under exchanging the first pair (or second pair) of 
indices; in our particular case, these symmetries are 
clear from the above formula 
for the curvature tensor, taken from [Totaro 04].  In this way, we obtain expressions 
of the required form for all the components of the curvature tensor.
Finally, we deduce that 
$$ \eqalign{ - 4 h \, R(\xi ,\eta ,\xi ,\eta ) 
&= - 4 h \sum _{{i<j,}\atop {p<q}} (\lambda _i \mu _j
- \lambda _j
\mu _i )  (\lambda _p \mu _q - \lambda _q \mu _p ) R_{ijpq} \cr & = 
S \, (\lambda _1 \mu _2 x_3 + 
\lambda _2 \mu _3 x_1  + \lambda _3 \mu _1 x_2 - \lambda _2 \mu _1 x_3 - 
\lambda _3 \mu_2 x_1 - \lambda _1 \mu _3 x_2 )^2 .\cr}$$

\noindent \bf 2. CONJECTURAL POSITIVITY OF $S$ FOR CERTAIN CUBICS ARISING IN GEOMETRY.

\rm In Section 5 of [Wilson 04], we were interested in the cubics which occur 
as intersection forms for 3-dimensional 
complete intersections in the product of three projective spaces.  We can 
however formalise this into a purely algebraic problem.  Suppose 
a ternary cubic is obtained as 
follows :  We choose positive integers $d_1 , d_2 , d_3$ and $r\ge 0$ such that 
$d_1 + d_2 + d_3 = r +3$, and set 
$$ P = (x_1 H_1 + x_2 H_2 + x_3 H_3)^3 \prod _{j=1} ^r ( a_j H_1 + 
b_j H_2 + c_j H_3 ),$$ 
with the $a_j$, $b_j$ and $c_j$ non-negative, and such that the cubic
$F(x_1 , x_2 , x_3)$, defined by taking the coefficient of the 
term in $H_1 ^{d_1} H_2 ^ {d_2} H_3
^{d_3}$ in the above formal product $P$, is non-degenerate. 
To relate this to the geometry, note that if the $a_j , b_j , c_j$ take 
non-negative integer values, then we may consider the complete intersection 
projective threefolds $X$ in $\P ^{d_1} \times \P^{d_2} \times \P ^{d_3}$ given 
by $r$ general trihomogeneous polynomials, with tridegrees 
$(a_j , b_j , c_j )$ for $j= 1, \ldots , r$.  The cubic we have defined above is then the
intersection form on the rank three sublattice of $H^2 (X, \Z )$ generated  by the pullbacks of
hyperplane classes from the three factors; by Lefshetz  
type arguments, this is usually the whole of $H^2
(X, \Z )$.

As in Section 1, we denote the coefficients of the ternary cubic $F$ by $a _{ijk}$,
where $i+j+k =3$.  These coefficients are themselves polynomials in the 
$a_j, b_j, c_j$, homogeneous of degree $r$ in each such set of variables.
We let  $S$ denote the $S$-invariant of $F$, and
$H$ the Hessian determinant of $F$. 

\proclaim  Conjecture 2.1.  Regarding $S$ as a polynomial in the $a_j, b_j, c_j$, 
every coefficient 
of this polynomial is non-negative.  \par
\proclaim Conjecture 2.2. Regarding $9 H^2 - 6^6 S F^2$ as a polynomial 
in the $a_j, b_j, c_j$ and $x_1 , x_2 , x_3$, every coefficient of this 
polynomial is also non-negative.\par

These conjectures imply their (weaker) geometric counterparts, in the case 
of $X$ being a complete intersection threefold in the product of 
three projective spaces, with second 
betti number three, and $F$ being its intersection form.   Here, we have 
taken specific non-negative integral 
values for the degrees $a_j, b_j, c_j$.  With the notation as in Section 0, these 
weaker conjectures may be interpreted, for $X$ as given, as saying that the curvature 
of the surface $\calK _1$ is bounded between $-9/4$ and zero.  The previous theoretical
and computational evidence for such conjectures to be true was outlined in Section 0
above.  Recall also that the first of these conjectures is equivalent to the statement 
that the semi-Riemannian Hessian metric on $\calK$ associated to $F$ has non-negative
Ricci curvature.

\proclaim Conjecture 2.3.  The intersection form of $X$ has non-negative $S$-invariant.\par
\proclaim Conjecture 2.4.  The polynomial $9 H^2 - 6^6 S F^2$ in $x_1, x_2, x_3$ takes 
non-negative values on the \Kahler cone of $X$, given by $x_1 > 0$, $x_2 > 0$, $x_3 > 0$.\par

Considered as 
a polynomial in the $a_j ,  b_j , c_j$, we have that $S$ is homogeneous of degree $4$ 
in any given set  $(a_j ,  b_j , c_j)$, and hence 
of total degree $4r = 4(d_1 + d_2 + d_3) -12$.  In fact, by inspection of the 
given formula for $S$, we see that $S$ is of degree $4d_1 - 4$ in the variables 
$(a_1 ,  \ldots , a_r)$, of  degree $4d_2 - 4$ in the variables 
$(b_1 ,  \ldots , b_r)$, and of  degree $4d_3 - 4$ in the variables 
$(c_1 ,  \ldots , c_r)$.

As explained above, the conjectures arose out of the theory developed in [Wilson 04];   
there is moreover now extensive computational evidence in their favour.  
In particular, Conjecture 2.1 has
been checked using MATHEMATICA for all $d_i \le 5$.  One can of 
course reduce to the case where all the $d_i$ equal some $d$, namely the maximum of the
$d_i$, by introducing $3d - d_1 - d_2 - d_3$ extra factors $(a_k H_1 + b_k H_2 + c_k
H_3)$ into the product $P$, and by considering the monomials in $S$, 
respectively 
$9H^2 - 6^6 S F^2$, which are of maximum 
possible degree 4 (respectively 6) in the appropriate variable $a_k$, 
$b_k$ or $c_k$ (and not therefore involving the other two).  For instance, for 
$d - d_1$ of the extra factors, the monomial considered should be of 
maximum degree 
in $a_k$ and not involve $b_k$ or $c_k$. 

To give 
a flavour of these calculations, I can report that in the case $d=3$ 
there are 209,520 non-zero terms in $S$, all with positive coefficients, 
and that a simple minded check of this took some two hours of computer time.
However, there are a very large number of symmetries, and taking such 
symmetries into account, the calculation was reduced to 
less than a couple of  
minutes.  For larger $d$ therefore, one should factor out by these symmetries.
For $d=4$ the conjecture was checked in a couple of  hours, 
and for $d=5$ in 
about four days.   The formula for $S$ given in Section 1 in terms of 
cofactors turns out to be slightly more efficient computationally than the formula 
in terms of the coefficients of the cubic.  The programs used by 
the author for these checks may
be found on his home page: http://www.dpmms.cam.ac.uk/$\sim$pmhw/
S{\_}invariant{\_}calculations.  The programs were run on a Sun V880 at the 
Max-Planck-Institut f\"ur Mathematik in Bonn, with 8 CPUs and 16GB of Main Memory 
theoretically available (although only a fraction of this would have been used).
The case $d=6$ seems to be beyond the range of standard 
computers.  The author has
not carried out as  many calculations on Conjecture 2.2, but it has been
verified for $d_1 = 3, d_2 = d_3 =2$, and there are strong 
 theoretical reasons in support of its geometric version (2.4), as outlined in 
Section 0.

  In this paper, we shall however concern ourselves mainly 
with the problem of Conjecture 2.1, that $S$ 
only has non-negative coefficients, 
and results closely related to this.

For the case $d=3$, one can obtain very precise information using MATHEMATICA 
about the coefficients.  The monomials appearing in any of the 25 terms in $S$ 
all appear in the expansion of $a_{111}^4$.  There are two types of monomial 
which appear in $a_{111}^4$ but not in $S$ (because the coefficients cancelling 
out) --- examples of these are 
$a_1 ^4 b_2 ^4 c_3 ^4 a_4 ^2 b_4 ^2 b_5 ^2 c_5 ^2 a_6 ^2 c_6 ^2$ and 
$a_1 ^4 b_2 ^4 c_3 ^4 a_4 ^2 b_4 c_4
a_5 b_5 ^2 c_5  a_6 b_6  c_6 ^2$.  If one considers the exponents as forming a 
$3 \times 6$ matrix, these monomials may be denoted rather more clearly as 
$$ \left[{\matrix{4&0&0&2&0&2 \cr 0&4&0&2&2&0 \cr 0&0&4&0&2&2 \cr}}\right] \quad 
\hbox{\rm and} \quad  
 \left[{\matrix{4&0&0&2&1&1 \cr 0&4&0&1&2&1 \cr 0&0&4&1&1&2 \cr}}\right] .$$ 
 Matrices 
differing from each other by permutations of the rows and/or columns 
are regarded as being of the same type.  There are then three types with 
coefficient 1 in $S$, represented by
matrices 
$$ \left[{\matrix{4&0&0&4&0&0 \cr 0&4&0&0&4&0 \cr 0&0&4&0&0&4 \cr}}\right] , \quad 
\left[{\matrix{4&0&0&3&0&1 \cr  0&4&0&1&3&0 \cr 0&0&4&0&1&3 \cr}}\right] , \quad
\left[{\matrix{4&0&0&0&2&2 \cr 0&4&1&1&2&0 \cr 0&0&3&3&0&2 \cr }}\right] .$$ 
A similar feature occurs for higher coefficients of there being rather 
a small number of types.  For instance, the largest coefficients which occur are 
$356$, $280$, $214$, $176$, $164$,
$128$, $106, \ldots$, all of which correspond to only one type.  
The highest coefficient $356$
corresponds to type 
$$ \left[{\matrix{2&1&1&2&1&1 \cr 1&2&1&1&2&1 \cr 1&1&2&1&1&2 \cr}}\right] .$$ 

It makes more sense however if we ignore all monomials 
containing fourth powers, on 
the grounds that these correspond to cases with the $d_i$ smaller.  With this
convention, the corresponding matrices do not have 4 in any entry.  In the case 
$d_1 = d_2 = d_3 =3$ as above, the smallest coefficients are then $4$, $6$ and $9$, 
corresponding (respectively) to matrices 
$$ \left[{\matrix{3&0&1&3&0&1
 \cr 1&3&0&1&3&0 \cr 0&1&3&0&1&3 \cr}}\right] , \quad 
\left[{\matrix{3&3&0&0&2&0 \cr  1&1&3&1&0&2 \cr 0&0&1&3&2&2 \cr}}\right] , \quad
\left[{\matrix{3&0&1&3&0&1 \cr 1&3&0&0&1&3 \cr 0&1&3&1&3&0 \cr }}\right] .$$ 
We shall also denote the first of these as $(d-1) \left[{\matrix{3&0&1
 \cr 1&3&0 \cr 0&1&3 \cr}}\right] , $ where $d=3$.

If we now move on to the cases $d_1 = d_2 = d_3 = d >3$, we may ask about 
the coefficients corresponding to $$(d-1) \left[{\matrix{3&0&1
 \cr 1&3&0 \cr 0&1&3 \cr}}\right] ; $$ for $d=4$, the coefficient may be 
calculated as 40, and for $d=5$ as 652.  For $d=4$ and $5$, 
a computer check verifies that this is the smallest non-zero 
coefficient (assuming no fourth powers) and the unique type of monomial 
corresponding to it, and one would conjecture that a similar statement is 
true for arbitrary values of $d>2$.  A formula for this coefficient for 
arbitrary $d$ will be produced in Section 4.

  Computer calculations suggest
 also a result that the cofactors  
$B_{pq}$ which appeared in Section 1 satisfy the condition 
that $B_{pq}$,  considered as a polynomial in 
the $a_i$, $b_j$, $c_k$ and $x_1 , x_2 , x_3$, 
 has only positive coefficients if $p \ne q$, and 
has only negative coefficients if $p=q$.    In the geometric situation of a 
three dimensional complete 
intersection in the product of three projective spaces, with the 
$(a_j , b_j , c_j )$ being assigned specific non-negative integral values, the
negativity of $B_{pp}$ corresponds to the Hodge index theorem on the 
surface cut out by $H_p =0$.  We shall 
prove these properties of the cofactors in the next section.

\noindent \bf 3. THE COFACTORS OF THE HESSIAN MATRIX.

\rm In this Section, we study further the cofactors $B_{pq}$ 
of the Hessian matrix of our ternary cubic $F$, where it will 
be more convenient here to denote the variables as $x_1, x_2, x_3$
rather than $x,y,z$.  
Recall that these cofactors were related to the $S$-invariant  by means 
of various expressions for $S$ described in Section 1; we shall return 
to this aspect in Section 4.
 In particular, 
for the special type of cubics we have studied in the last two sections, 
the $B_{pq}$ may be considered as polynomials in the $a_j$, $b_j$, $c_j$ and 
$x_1$, $x_2$, $x_3$.  In this Section, 
we confirm the expectations, mentioned in Section 2, 
  concerning the signs
of their coefficients; this in turn will show that the Hessian determinant $H$
only has positive coefficients (4.2).  This latter fact might be expected because of
the Hodge Index Theorem, which implies the weaker statement that $H$ takes 
non-negative values for non-negative values of its variables.

\proclaim Theorem 3.1.  The polynomials $B_{pp}$ only have negative coefficients, 
and the polynomials $B_{pq}$ for $p\ne q$ only have positive 
coefficients.\par

\noindent\it Proof.\rm\quad  For the first part, we may consider $B_{33} 
= f_{11} f_{22} - f_{12}^2$.  
For a general cubic $f$, we have 
$$ {1\over {36}} B_{33} = (a_{300} x_1 + a_{210} x_2 + a_{201} x_3 )
( a_{120} x_1 + a_{030} x_2 + a_{021} x_3 ) - 
(a_{210} x_1 + a_{120} x_2 + a_{111} x_3 )^2 .$$   The fact
that, in our particular case,  this polynomial is non-positive for all non-positive 
values of the variables follows from the Hodge index theorem again.  We 
however prove the more precise result that the coefficients are all negative.

Let us consider for instance the term in $x_1 x_2$;  we prove that its 
coefficient $$a_{300} a_{030} - a_{210} a_{120},$$ 
 considered as a polynomial 
in the $a_i$, $b_j$, $c_k$, has only negative coefficients.  Without loss 
of generality, we can assume that $d_1 = d_2 = d_3 = d$, and we set $s = d-1$.
Then the polynomial in question is of degree 2 in each set of variables 
$(a_j , b_j , c_j )$, and is of degree $2s-1$ in the $a_i$, degree 
$2s-1$ in the $b_j$, and degree $2s +2$ in the $c_k$.  On the other hand, $a_{300}$ 
(respectively, $a_{210}$) is of degree $s-2$ (respectively, $s-1$) in the $a_i$, 
degree $s+1$ (respectively, $s$) in the $b_j$, and degree $s+1$ 
(respectively, $s+1$) in the $c_k$, with analogous statements for $a_{030}$ 
and $a_{120}$.

Let us now consider a monomial of the appropriate degrees in the 
$(a_j ,b_j , c_j )$, and ask about its coefficient as a term in 
$a_{300} a_{030} - a_{210} a_{120}$.  We suppose that the monomial in 
question consists of $p_1$, respectively $p_2$, $p_3$, 
occurrences (for various $j$) of $a_j ^2$, respectively  
$b_j ^2$, $c_j ^2$, and  $\tilde u$, respectively $\tilde v$, $\tilde w$,  
occurrences of $a_j b_j$, 
respectively $a_j c_j$, $b_j c_j$.  As in 
Proposition 5, we shall see that only the mixed cases will be of relevance. 
Note that $2p_1 + \tilde u + \tilde v = 2s -1$, $2 p_2 + \tilde u + 
\tilde w = 2s -1$ and  $2p_3 + \tilde v + \tilde w = 2s +2$.

The coefficient of the monomial in $a_{300} a_{030}$ 
is given by counting the number of ways of
expressing it as a monomial in $a_{300}$ times a monomial in $a_{030}$, 
and similarly for its coefficient in $a_{210} a_{120}$.  To obtain the 
first factor in the former case,  
involves choosing $s-p_1 -2 = \ha (\tilde u+ \tilde v-3)$ of the 
$a_j b_j$ and  $a_j c_j$ appearing for which we choose the $a_j$, $s-p_2 +1 = \ha (
\tilde u+ \tilde w +3)$ of the 
$a_j b_j$ and  $b_j c_j$ for which we choose the $b_j$, and 
$s - p_3 +1 = \ha (\tilde v+ \tilde w)$ of the $a_j c_j$ and $b_j c_j$ for which we 
choose the $c_j$.  Note here the necessary parity condition that either 
$\tilde u$ is odd and $\tilde v,\tilde w$ are even, 
or the other way round.  We shall deal with 
the first case; the other case follows similarly.

We set $\tilde u = 2 u +1$, $\tilde v = 2 v$ and $\tilde w = 2 w$.  The 
possible factorizations are then given by choosing $k$ of the $2u+1$ 
occurrences of $a_j b_j$ for which we choose the $a_j$, choosing 
$u+v -1 -k$ occurrences of $a_j c_j$ for which we choose the $a_j$, and 
finally $w-u +1 +k$ occurrences of the $b_j c_j$ for which we choose the 
$b_j$, the rest then being determined.   Thus the number of ways of doing 
this, and hence the coefficient of the monomial in $a_{300} a_{030}$, is
$$ \sum _{k=0} ^{2u+1} {{2u +1} \choose k} {{2v}\choose {v+u -(k+1)}} {{2w} 
\choose { w + u - (k+1)}}.$$   Similarly, the coefficient of the 
monomial in $a_{210} a_{120}$ is seen to be 
$$ \sum _{k=0} ^{2u+1} {{2u +1} \choose k} {{2v}\choose {v+u - k}} {{2w} 
\choose { w + u - k }}.$$  Thus we need to verify the negativity of  
$$  \sum _{k=0} ^{2u+1} {{2u +1} \choose k} \left( {{2v}\choose {v+u -(k+1)}} {{2w} 
\choose { w + u - (k+1)}} - {{2v}\choose {v+u - k}}{{2w} 
\choose { w + u - k }} \right).$$  

This sum may however be rearranged as 
$$ \eqalign{& - {{2v}\choose {v+u}} {{2w}\choose {w+u}} + 
{{2v}\choose {v+u+2}} {{2w} \choose {w+u+2}} 
\cr &
- \sum _{k=1} ^{2u+1} {{2v}\choose {v+u -k}} {{2w}\choose {w+u -k}} 
\left( {{2u+1}\choose {k}} - {{2u+1}\choose {k-1}} \right) .\cr }$$ 
The first line of this rearranged sum is now clearly non-positive.  
In the summation, the term 
$\left( {{2u+1}\choose {k}} - {{2u+1}\choose {k-1}} \right)$ is 
antisymmetric about $u+1$, and in fact equals 
$$   {{2(u+1 -k)}\over {2u+2}} {{2u+2}\choose k}.$$  
If we therefore pair these antisymmetric terms, and 
use the fact that  for $j>0$, we have 
$$ {{2v} \choose {v+j-1}} \ge {{2v} \choose {v+j+1}}, \quad 
 {{2w} \choose {w+j-1}} \ge {{2w} \choose {w+j+1}},$$
the claimed inequality follows.

For the term in $x_1 x_3$, we need to show that the polynomial 
$a_{300} a_{021} + a_{201} a_{120} - 2 a _{210} a_{111}$ only 
has negative terms.  In fact, we prove this for the two polynomials 
$a_{300} a_{021} -  a _{210} a_{111}$ and 
$a_{201} a_{120} -  a _{210} a_{111}$.  Let us consider a particular 
monomial appearing in these polynomials; with the notation as above, 
the parities on $\tilde u$, $\tilde v$
and $\tilde w$ will differ from before.  Since $2p_1 + \tilde u + \tilde v = 2s-1$, 
$2p_2 + \tilde u + \tilde w = 2s$ and $ 2 p_3 + \tilde v + \tilde w = 2s+1$, 
we have either $\tilde v$ odd and $\tilde u$, $\tilde w$ even, 
or the other way round.  Considering for instance the case $\tilde v = 2v + 1$, 
$\tilde u = 2u$ and $\tilde w = 2w$, we can run through a similar argument 
to that given above, and find that the coefficient 
of the given monomial in $a_{300} a_{021} -  a _{210} a_{111}$ is 
$$  \sum _{k=0} ^{2v+1} {{2v +1} \choose k} {{2w}\choose {w+v-k}}
\left( {{2u}\choose {u +v  -(k+1)}}  - {{2u}\choose {u +v - k}}\right).$$
We now observe that the bracket in this summation is antisymmetric 
about $k = v - \ha$, and then pairing off terms proves the result in an 
analogous 
way to before.  Similarly, the coefficient of the given monomial in 
$a_{201} a_{120} -  a _{210} a_{111}$ is 
$$  \sum _{k=0} ^{2v+1} {{2v +1} \choose k} {{2u}\choose {u+v-k}}
\left( {{2w}\choose {w +v  -(k+1)}}  - {{2w}\choose {w +v - k}}\right),$$ 
and the same argument goes through, switching the roles of $u$ and $w$.

For the term in $x_1 ^2$, we need to show that the polynomial 
$a_{300} a_{120} - a_{210}^2$ only has negative terms.  For a monomial 
to appear in this polynomial, we have yet another parity condition, namely
that $\tilde u$, $\tilde v$ and $\tilde w$ are all even, or are all odd.  
The reader is left to 
check the negativity.  By symmetry, the only other term we need to consider is 
that in $x_3 ^2$; here we need that the polynomial 
$a_{201} a_{021} - a_{111}^2$ only has negative terms.  
The parity condition here is the same as for $x_1
^2$, and the reader is left to verify the details of the negativity.

We now need to consider the cofactors $B_{pq}$ with $p\ne q$.  We shall only 
explicitly verify the $x_3^2$ terms here, and leave the others to the reader.
Note in passing that in the formula for $6^4 S x_3^2 $ from Section 1,  
we may consider instead the identity given simply by the terms in $x_3^2$, and 
so it will be the $x_3 ^2$ terms in the above cofactors which will occur in the 
algorithm we describe in Section 4.
 We check these terms for $B_{12}$ and $B_{13}$, 
the rest then following from considerations of symmetry.  Let us 
start with ${1\over {36}} B_{12}$, which is 
$$\eqalign{{1\over {36}} (f_{13} f_{23} - f_{12} f_{33}) =  & \ 
 ((a_{201} x_1 + a_{111} x_2 + a_{102} x_3) 
(a_{111} x_1 + a_{021} x_2 + a_{012} x_3) \cr & - 
(a_{210} x_1 + a_{120} x_2 + a_{111} x_3)(a_{102} x_1 + a_{012} x_2 
+ a_{003} x_3 )),\cr }$$ whose term in $x_3^2$ is 
$$ a_{102} a_{012} - a_{111} a_{003}.$$  For ${1\over {36}} B_{13} = 
{1\over {36}} (f_{12} f_{23} - f_{13} f_{22})$, we have instead the 
polynomial $$ a_{111} a_{012} - a_{102} a_{021}.$$  The latter we already know
has only positive terms from our calculations on the $x_1 x_3$ term for 
$B_{33}$, where we saw that the polynomial $a_{201} a_{120} - a_{210} a_{111}$ 
only had negative terms (simply switch the first and last indices).  
For a given monomial to appear in the first polynomial, 
we need parities that $\tilde u$ is odd and
$\tilde v$, $\tilde w$ even, or the other way round.  
For the monomial to appear in the second
polynomial, we need parities that 
$\tilde v$ is odd and $\tilde u$, $\tilde w$ even, or the other way round.

For the former, namely $ a_{102} a_{012} - a_{111} a_{003}$, the by now 
familiar calculation shows that the coefficient of our monomial, say in the case 
$\tilde u = 2u+1$ odd and $\tilde v = 2v$, $\tilde w = 2w$ even, is the sum 
$$  \sum _{k=0} ^{2u+1} {{2u +1} \choose k} {{2v}\choose {v+u-k}}
\left( {{2w}\choose {w + u  - k}}  - {{2w}\choose {w +u - (k-1)}}\right).$$ 
The bracketed term is now antisymmetric about $k = u+\ha $, and pairing the 
terms again, we see that the sum is positive.
\bigskip
\proclaim Theorem 3.2.  For the cubics under consideration, the Hessian determinant $H$ 
is a polynomial in the $a_j$, $b_j$, $c_j$ and 
$x_1$, $x_2$, $x_3$, all of whose coefficients are positive.\par

\noindent \it Proof.\rm\quad
Recall that, for any $n \times n$ matrix $A$ with $n>2$, 
we have $\Adj ( \Adj A ) = 
\det (A)\, A$.  Applying this, with $A = (F_{ij})$, we deduce that
$$ F_{12} H  = - B_{33} B_{12} + B_{23} B_{13}.$$
Theorem 3.1 then implies that $F_{12} H$, a polynomial in the $a_j$, $b_j$, $c_j$ and 
$x_1, x_2, x_3$, only has positive coefficients, 
where we may without loss of generality assume that $F_{12}$ is non-trivial.

An easy argument shows however that if $f,g$ are polynomials in a 
finite set of variables, with $f$ non-trivial, 
such that $f$ and $fg$ only have positive 
coefficients, then the same is true for $g$.  To see this, choose an
order for the variables, and then order the monomials lexicographically.  
Now pick the largest monomial in $f$, and the largest monomial (if it 
exists) whose coefficient in $g$ is negative; the product of these terms 
would yield a monomial in $fg$ with negative coefficient.  Applying this, 
since $F_{12}$ only has positive coefficients, we deduce that the same
holds for $H$.\medskip

\noindent \bf Remark 3.3. \quad\rm  Once we know that $H$ only has positive
coefficients, then the above argument shows that the same is true for all 
entries of $\Adj (B)$, for instance $B_{11} B_{22} - B_{12} ^2 = F_{33} H$.  
\bigskip
\eject

\noindent \bf 4. MORE COMBINATORICS OF THE $S$-INVARIANT.

\rm The fact that we have explicit formulae for the coefficients in both 
the polynomials 
$\partial ^2 B_{33}/\partial x_p \partial x_q$ and $\partial ^2 B_{pq} /
\partial  x_3 ^2$, provides an explicit recipe for calculating the 
coefficients in $S$ directly.  From 
Remark 1.2, we know that 
$$  \qu \sum _{p, q} (\partial ^2 B_{pq} / \partial x_3 ^2 )
(\partial ^2 B_{33}/ \partial x_p \partial x_q) = 
6^4 S .$$  
The tridegrees of the terms in $\partial ^2 B_{33}/ \partial x_p \partial x_q$ are 
$$ \pmatrix{ (2s-2 , 2s, 2s+2) & (2s-1 , 2s-1, 2s+2) & (2s-1 , 2s , 2s+1) \cr
(2s-1,2s-1,2s+2) & (2s, 2s-2, 2s+2) & (2s , 2s-1, 2s+1) \cr
(2s-1,2s, 2s+1) & (2s, 2s-1, 2s +1) & (2s,2s,2s) \cr }, $$ and those of 
$\partial ^2 B_{pq}/\partial x_3 ^2$ the complementary degrees 
with respect to $4s$; for instance the tridegree of 
 $\partial ^2 B_{12}/ \partial x_3 ^2$ is $(2s+1, 2s+1, 2s-2)$.
  The 
recipe for calculating the coefficient of a given allowable monomial $M$ 
in now clear.  Consider all factorisations $ M = M_1 M_2 $ of $M$, 
where the $M_i$ are quadratic in each set of variables $(a_j , b_j , c_j )$, 
and where 
$M_2$ has one of the tridegrees listed above for  
$\partial ^2 B_{33}/ \partial x_p \partial x_q$, 
with $M_1$ having the complementary tridegree.
The $M_i$ give rise to numbers $\tilde u_i , \tilde v_i , \tilde w_i$, 
where $i = 1,2$, from which we have an explicit expression for the 
coefficient of $M_i$ in the relevant entry of the matrix in question.  
Adding the products of these two coefficients as we range over the
 factorizations gives 
us the coefficient of $M$ in $S$.

We illustrate this with the following example;
we consider the case $s = 3 t$, and 
so $d = 3t +1$, and $M$ a monomial with matrix of exponents having $4t$ 
columns of the form $\pmatrix{3\cr 1 \cr 0\cr }$, $4t$ columns of the form 
$\pmatrix{0\cr 1 \cr 3 }$, and $t$ columns of the form 
$\pmatrix{0\cr 4 \cr 0}$.    Note that for all factorizations $M = M_1 M_2 $, 
we have $\tilde v_1 = 0 = \tilde v_2$.   A factorization is 
determined by specifying for how many of the $a_j^3 b_j$ one takes $a_j ^2$ in $M_1$, and 
for how many of the $b_j c_j^3$ one takes $c_j^2$; 
if these numbers are denoted by $k,l$ respectively, 
then $\tilde u_1 = 4t -k$, $\tilde w_1 = 4t -l$, 
$\tilde u_2 = k$, $\tilde w_2 =l$.  Note that $M_2$ then has tridegree 
$(8t -k, 2t +k+l, 8t-l)$.  Thus the only pairs $(k,l)$ of relevance will
be $(2t, 2t)$, $(2t, 2t-1)$, $(2t, 2t-2)$, $(2t+1, 2t -1)$, $(2t+1, 2t -2)$ 
and $(2t+2, 2t-2)$.  We consider each pair in turn; the fact that $\tilde v 
= 0$ simplifies the algebra considerably. The case $(2t ,2t)$ corresponds to 
the $x_3 ^2$ term in $B_{33}$; the coefficient 
of the monomial in $B_{33}/36$ is checked to simplify to
$$ {{2t}\choose {t}} \left( {{2t}\choose {t-1}} - {{2t} \choose {t}} \right).$$
The case $(2t ,2t-1)$ corresponds to the $x_2 x_3$ term in $B_{33}$; 
the relevant
coefficient is $$2{{2t -1}\choose t} \left( {{2t} \choose {t-1}} - 
{{2t} \choose {t}}\right) .$$  The case $(2t, 2t-2)$ corresponds to the $x_2^2$ term in
$B_{33}$; the relevant 
coefficient is $${{2t - 2}\choose {t-1}} \left( {{2t} \choose {t-1}} - 
{{2t} \choose {t}}\right) .$$  The case 
$(2t+1,2t-1)$ corresponds to the $x_1 x_3$ term; the relevant coefficient is 
 $${{2t - 1}\choose {t}} \left( {{2t +1} \choose {t-1}} - {{2t +1} \choose {t}}\right) .$$   
The case $(2t+1, 2t-2)$ corresponds to the $x_1 x_2$ term;
the relevant coefficient is  $${{2t - 2}\choose {t-1}} 
\left( {{2t +1} \choose {t-1}} - {{2t +1} \choose {t}}\right) .$$   
The case $(2t+2, 2t-2)$ corresponds to the 
$x_1 ^2$ term; the relevant coefficient is  
$${{2t - 2}\choose {t-1}} \left( {{2t +2} \choose {t}} - {{2t +2} \choose {t +1}} \right) .$$

Now we need the corresponding $x_3 ^2$ terms in $B_{pq}$.  We already know
that $(k,l) = (2t,2t)$ corresponds to the $x_3^2$ term in $B_{33}$ with 
coefficient of the monomial in $B_{33}/36$ being 
$${{2t}\choose {t}} \left( {{2t} \choose {t -1}} - {{2t} \choose {t}}\right) .$$  
We check that $(2t,2t-1)$ corresponds to the $x_3^2$ term 
in $B_{23} = f_{12}f_{13} - f_{11} f_{23}$, namely 
$36 (a_{111} a_{102} - a_{201} a_{012} )$, and that the coefficient required is 
$${{2t +1}\choose {t}} \left( {{2t} \choose {t}} - {{2t} \choose {t -1}}\right)  ;$$  
the pair $(2t,2t-2)$ corresponds to the $x_3^2$ term of 
$B_{22} = f_{11} f_{33} - f_{13}^2$, namely $36(
a_{201} a_{003} - a_{102}^2 )$, and the coefficient is
$${{2t}\choose {t +1}}{{2t +2} \choose {t+2}} - {{2t} \choose {t}}{{2t+2} 
\choose {t+1}} ;$$ the pair 
 $(2t+1,2t-1)$ corresponds to the $x_3^2$ term of 
$B_{13} = f_{12} f_{23} - f_{13} f_{22}$, namely $36( a_{111} a_{012} - a_{102} a_{021})$, and the
coefficient is
$${{2t-1}\choose {t -1}} {{2t +1} \choose {t+1}} - {{2t-1} \choose {t-1}}
{{2t+1} \choose {t}} = 0;$$ the pair 
 $(2t+1,2t-2)$ corresponds to the $x_3^2$ term of 
$B_{12} = f_{13} f_{23} - f_{12} f_{33}$, namely $36( a_{102} a_{012} - a_{111} a_{003})$, and the
coefficient is
$${{2t-1}\choose {t -1}} \left( {{2t +2} \choose {t+1}} - 
{{2t+2} \choose {t+2}}\right) ;$$  finally  $(2t+2,2t-2)$ 
corresponds to the $x_3^2$ term of $B_{11} = f_{22} f_{33} - f_{23}^2$, namely $36( a_{021} a_{003} -
a_{012}^2 )$, and the coefficient is
$${{2t-2}\choose {t -1}} \left( {{2t +2} \choose {t}} - 
{{2t+2} \choose {t+1}}\right) .$$  

We now have all the information we need to 
calculate $S$ from the  
formula given at the start of the Section, where of course for a given $(k,l)$ we shall
need to weight the contribution by $ {{4t}\choose k} {{4t}\choose l}$.  
Putting all this together, we get a formula for the 
relevant coefficient of the $S$-invariant as a 
function of $t$.  
With the aid of MATHEMATICA, one can then simplify 
the formula to the surprisingly simple form 
$$ (4t)!^2 \left( {t! \over { (t-1)! (t+1)! }} - {1 \over {t!^2}} \right) ^4.$$
In particular, one notes that it is positive.
Evaluating this formula for $t$ taking values $0,1,2,3,4, \ldots $, one 
obtains values for the coefficient of the monomial in the $S$-invariant to be
$1, 36, 78400$, $533610000$, $6363107150400, \ldots $.  The first two 
of these values coincide with previously calculated numbers (using a simple-minded
method).

The author has checked positivity of the coefficient for other cases of a similarly 
general type.  Apart from the computer calculations described in Section 2, 
perhaps the most
telling evidence however 
for the positivity of all the coefficients is provided by calculating what
was conjectured in Section 2 to be the \it smallest \rm coefficient.

\noindent \bf Example 4.1.\rm\quad We consider therefore 
the case where $s= d-1$ and the 
monomial $M$ has a matrix of 
exponents  
 $$(d-1) \left[{\matrix{3&0&1
 \cr 1&3&0 \cr 0&1&3 \cr}}\right] .$$  The 
coefficient was calculated for $d\le 5$ in 
Section 2, and the monomial was conjectured to have the smallest 
coefficient (assuming no 
fourth powers) for any given value of $s$.   For this reason, it is an obvious 
crucial case in which to verify our main conjecture.  
In a factorization 
$M = M_1 M_2$, we suppose that for precisely $k$ of the $a_j ^3 b_j$ we have 
taken $a_j ^2$ in $M_1$, for precisely $l$ of the $a_j c_j ^3$ we have taken 
$c_j ^2$ in $M_1$, and for precisely $m$ of the $b_j ^3 c_j$ we have taken $b_j ^2$.  
Consideration of tridegrees shows that the only pairs $(k-l, m-l)$ of 
relevance are $(0,0)$, $(0,1)$, $(0,2)$, $(1,1)$, $(1,2)$ and $(2,2)$.  For a 
given choice of $(k,l,m)$, the corresponding 
triple $(\tilde u , \tilde v , \tilde w )$ associated with $M_2$ is just $(k,l,m)$.  
Because the $\tilde v$ is no longer 
zero in general, the formula for the coefficient (as a function of $s$) 
that we obtain involves triple summations.  The 
rather complicated formula (occupying a page) which results may be found in an 
Appendix to this paper.
Although MATHEMATICA does not reduce this formula to any simple form, it is  
nevertheless an explicit formula, which has been checked to give positive 
values for $s \le 501$.  The proof of positivity for general $s$ 
presumably follows by suitably rearranging the 
sums which occur in the formula.   
The values for $s= 1,2,3,4,5,6,7,8$ are 
respectively $1$, $4$, $40$, $652$,
$13174$, $308464$, $8158021$ and $23830660$; the first four of these 
correspond to the previously calculated
values.   The  fact
that the numbers generated tend to have large prime factors 
(for instance $8158021$ is prime) suggests
that there is no simple form of the formula.  We should also comment that the cofactor
formula for $S$ that we are using expresses the coefficient of the given monomial as
the sum of six terms.  By taking for instance $s=4$ in this example and evaluating these 
terms, each of the terms has modulus greater than the sum of the terms;  so although
the sum is positive, it does involve significant cancellations.

A proof of the positivity of the coefficient for the case of a general 
monomial still seems some way off, at least using the recipe given above.
I restrict myself to the comment that the formulae we derived for 
the coefficients of monomials in the cofactors can all be expressed as the 
difference between two reasonably simple hypergeometric series of the form 
${}_3 F_2$ --- in some of the special cases worked out, they were the 
difference of even simpler terms.  The theory of hypergeometric series 
may therefore feature in a proof of the conjectures and in possible alternative proofs 
of the results from Section 3.

It might be observed that there are other relatively simple formulae which yield $S$,
apart from those in Section 1.  By a similar method of proof to Lemma 1.1, one can for
instance show  that $$ \ha \sum B_{ij}\, \partial ^2 H /\partial x_i \partial x_j = 6^5
S F .$$  If one could prove positivity of the coefficients for this polynomial, then the
desired result would follow for $S$.  With the methods described above however, the
expression for $S$ that we have used is simpler to analyse than this one.

\noindent \bf 5. CONCLUDING REMARKS.

\rm  We summarised in Section 0 the theoretical evidence for the geometric conjectures
(2.3) and (2.4), that for complete intersection threefolds in the product of three
projective spaces (assuming $b_2 =3, h^{2,0} =0$), 
the curvature of the surface $\calK _1$ is bounded between $-9/4$ and 0.  These
conjectures were set in the more general context of threefolds admitting certain
specific types of degeneration, and for $b_2 \ge 3$ can be rephrased in an illuminating
way in terms of Ricci curvatures.  In the Calabi--Yau case, there was further evidence
via mirror symmetry from known results on the Weil--Petersson metric on the complex
moduli space of the mirror.

Even if we knew however that Conjectures 2.3 and 2.4 held, it is unclear whether this
would help in a proof of their combinatorial versions (2.1) and (2.2).  An illustration 
of this is that the geometric version of Theorem 3.2 follows from the Hodge index
theorem, but this does not seem to help in a proof of the combinatorial result, or in 
a proof of Theorem 3.1.  If one could produce a proof of (3.1) which depended 
less on explicit combinatorial manipulations than the proof given here,
I believe that this might suggest alternative approaches to proofs of (2.1) and (2.2).

The experimental evidence for Conjecture 2.1 is I believe very strong indeed.  Not only
has it been checked in all cases up and including $d=5$, this involving a 
prodigious amount of calculation, but it has also been checked in the case of the
predicted minimum coefficient up to enormous values of $d$.  The computational 
evidence for Conjecture 2.2 is not as extensive, although still strong, but there is  
more theoretical evidence in the geometric case
for the precise value of the upper bound.  Should one want
further experimental evidence for (2.2), it should be feasible  to extend
the previous computations at least to include all cases up to and including $d=3$.

\noindent \bf ACKNOWLEDGEMENT

\rm This paper was written during a stay at the Max-Planck-Institut in 
Bonn, which 
he would like to thank for both financial and computational support.  

\noindent\bf REFERENCES 

\rm
\noindent [Aronhold 58] \quad S. Aronhold,  Theorie der homogenen Functionen 
dritten Grades von drei \break Ver\"anderlichen, \it  J. reine angew. Math. 
\bf  55 \rm (1858) 
97-191.

\noindent [Lu 01] \quad Zhiqin Lu, On the Hodge metric of the universal deformation 
space of Calabi--Yau threefolds, \it J. Geom. Analysis \bf 11 \rm (2001) 103-118.

\noindent [Lu \& Sun 04] \quad Zhiqin Lu \& Xiaofeng Sun, Weil--Petersson geometry on
moduli space of 
polarized Calabi--Yau manifolds, \it J. Inst. Math. Jussieu \bf 3 \rm (2004) 185-229.

\noindent [O'Neill 83] \quad Barrett O'Neill, \it Semi-Riemannian Geometry \rm 
(Academic Press, San Diego, 1983).

\noindent [Sturmfels 93] \quad B. Sturmfels, \it Algorithms in invariant theory \rm 
(Springer, Wien, 1993).

\noindent [Totaro 04] \quad B. Totaro,   The curvature of a Hessian metric, 
\it  Int. J. Math. \bf 15 \rm (2004) 369-391.

\noindent [Wilson 92] \quad P.M.H. Wilson,  The K\"ahler cone on Calabi--Yau threefolds, 
\it Invent. math. \bf 107 \rm (1992) 561-583; Erratum: \it Invent. math. \bf 114 \rm
(1993) 231-232.

\noindent [Wilson 04] \quad P.M.H. Wilson,  Sectional curvatures of K\"ahler moduli, 
\it Math. Ann. \bf 330 \rm (2004) 631-664.   

\bigskip 

\noindent \bf Appendix.\rm

 The formula for the coefficient $A$ of $M$ in $S$, with $M$ as in Example 4.1, 
is given as $A = A_1 +  A_2 + A_3 + A_4 + A_5 + A_6$, where the $A_i$ are  
defined (as functions of $s$) as follows:
$$ A_1 = \sum _{l=0} ^s \sum _{j=0} ^l \sum _{i=0} ^{s-l} {s \choose l}^3 
{l\choose j}^2 \left( {l\choose j+1} - {l\choose j} \right) {s-l\choose i}^2 
\left( {s-l\choose i+1} -
{s-l\choose i} \right) .
$$
$$ \eqalign{ A_2 = \sum _{l=0} ^{s-1}  \sum _{j=0} ^l  
{s \choose l}^2 & {s \choose l+1} {l\choose j} \left( {l\choose j+1} 
{l+1\choose j+2} + {l\choose j+1} {l+1\choose j+1} - 2 {l\choose j} 
{l+1\choose j+1}
\right)\cr & \sum _{i=0} ^{s-l}{s-l \choose i} \left( {s-l\choose i}  
{s-l-1 \choose i-1} - {s-l \choose i+1} {s-l -1\choose i}\right) .\cr}
$$
$$ \eqalign{ A_3 = \sum _{l=0} ^{s-2}  \sum _{j=0} ^l 
{s \choose l}^2 & {s \choose l+2} 
{l\choose j} \left( {l\choose j+1} {l+2\choose j+2} - {l\choose j} 
{l+2\choose j+1}\right) \cr & 
\sum _{i=0} ^{s-l}{s-l \choose i}  {s-l-2\choose i-1}  
\left( {s-l \choose i+1} - {s-l \choose i}\right) .\cr}$$
$$ \eqalign{ A_4 = \sum _{l=0} ^{s-1}  \sum _{j=0} ^{l +1}
{s \choose l+1}^2 & {s \choose l} 
 \left( {l\choose j+1} {l+1\choose j+1} + {l\choose j} 
{l+1\choose j+1} - 2 {l+1\choose j} {l \choose j-1 }\right) 
\cr & {l+1 \choose j}\sum _{i=0}
^{s-l-1}{s-l-1 \choose i}  {s-l\choose i+1}  \left( {s-l-1 \choose i} - 
{s-l-1 \choose i+1}\right) .\cr}$$
$$ \eqalign{ A_5 = \sum _{l=0} ^{s-2}  \sum _{j=0} ^{l}
& {s \choose l} {s \choose l+1} {s \choose l+2}
 {l +1\choose j}  \left( {l\choose j+1} {l+2\choose j+2} -  
{l\choose j} {l+2  \choose j+1 }\right) \cr & 
\sum _{i=0} ^{s-l-1} {s-l-1 \choose i}  {s-l\choose i+1} 
\left( {s-l-2 \choose i} - {s-l-2 \choose i-1}\right) .\cr}$$
$$ \eqalign{ A_6 = \sum _{l=0} ^{s-2}  \sum _{j=1} ^{l+2}
{s \choose l+2}^2 & {s \choose l}
 {l +2\choose j}  \left( {l\choose j-2} {l+2\choose j-1} -  
{l\choose j-1} {l+2  \choose j }\right) \cr & \sum _{i=0} ^{s-l-2} 
{s-l-2 \choose i}  {s-l\choose i+1} 
\left( {s-l-2 \choose i-1} - {s-l-2 \choose i}\right) .\cr}$$

If we take the formula for $S$ in terms 
of cofactors, as used in Section 4,
but write it as a sum over $p \le q$, 
these numbers represent the 
coefficients of $M$ in the terms 
with $(p,q) = (3,3)$, $(2,3)$, $(2,2)$, $(1,3)$, $(1,2)$ and $(1,1)$,  
respectively. If we take as an example $s=4$ in the given formulae, the above
numbers are $A_1 = 5804$, $A_2 = -3048$, $A_3 = 
2352$, $A_4 = -4552$, $A_5 = -2256$, $A_6 = 2352$ and $A = 652$.  In 
fact, for the monomial $M$ of this example, we have $A_3 = A_6$ for all $s$;
this latter identity may be seen by writing $A_3$ in terms of $l' = s-2-l$, 
rearranging the sums over $i$ and $j$,  
and then comparing with the formula for $A_6$.

\end